\documentclass{jsiamletters}
\group{
Algorithms for Matrix / Eigenvalue Problems and their Application
}
\newcommand{\newblock}{\hskip 1mm}

\usepackage{color}
\newcommand{\rr}[1]{#1}
\newcommand{\rrr}[1]{#1}

\affiliation{Department of Mechanical and Physical Engineering, Tottori University}{Tottori,  680-8552, Japan}
\affiliation{Information Initiative Center, Hokkaido University}{Sapporo, 060-0811, Japan}
\affiliation{Department of Communication Engineering and Informatics, The University of Electro-Communications}{Chofu, 182-8585, Japan}

\authorinfo{Hisashi Kohashi}{1}{d20t1003m@edu.tottori-u.ac.jp}
\authorinfo{Harumichi Iwamoto}{1}{b18t1015b@edu.tottori-u.ac.jp}
\authorinfo{Takeshi Fukaya}{2}{fukaya@iic.hokudai.ac.jp}
\authorinfo{Yusaku Yamamoto}{3}{yusaku.yamamoto@uec.ac.jp}
\authorinfo{Takeo Hoshi}{1*}{d20t1003m@edu.tottori-u.ac.jp}
\email{hoshi@tottori-u.ac.jp}
\title{Performance prediction of massively parallel computation by Bayesian inference}
\abstract{
A performance prediction method for massively parallel computation is proposed. The method is based on  performance modeling and Bayesian inference to predict  elapsed time $T$ as a function of the number of used nodes $P$ ($T=T(P)$). 
The focus is on extrapolation for larger values of $P$ from the perspective of application researchers. The proposed method has several improvements over the method developed in a previous paper, and application to real-symmetric generalized eigenvalue problem  shows promising prediction results. The method is generalizable and applicable to many other computations. }
\begin{document}

\maketitle

\section{Introduction}

A number of parallel algorithms for modern supercomputers
have been developed in recent decades. 
However,it remains difficult
for the general application researcher
to choose the proper algorithm and/or set the values of tuning parameters
that influence computational performance.
One practical remedy is 
to predict the elapsed time \rr{$T$} 
as \rr{a} function of the number of processor nodes $P$ ($T=T(P)$).
In particular, 
application researchers \rr{are interested in} extrapolation for larger numbers of nodes. 

The present paper proposes a performance prediction method
featuring improvements over the method suggested 
in our previous paper ~\cite{EigenKernel}. 
In Ref.~\cite{EigenKernel}, 
we presented several generic performance models 
for the real-symmetric generalized eigenvalue problem in which
the performance model function $T=T(P)$ 
is based on the applied algorithm and architecture
and contains several fitting parameters.  
The parameters are estimated by Bayesian inference, 
whereby the elapsed time is predicted or extrapolated from the teacher data or the existing benchmark data. 
The method in the present paper improves that proposed in  Ref. \cite{EigenKernel}  in two respects:  
(i) it improves the performance model and the cost function for the Bayesian inference and
(ii) it allows for a systematic comparison among different models and different teacher data set.
Here, we use the same performance data as used in Ref. \cite{EigenKernel},
since it shows the typical performance behavior of parallel computation 
and will enable the reader to compare the present results with those in Ref. \cite{EigenKernel}.

The method and results are described in Sections  \ref{SEC-METHOD} and  \ref{SEC-RESULTS}, respectively.  
A summary is provided in Section \ref{SEC-SUMMARY}.  

\section{Method \label{SEC-METHOD}}

\subsection{Performance models}

Performance models were investigated in papers
such as Ref.~\cite{Dackland, Grbovic07,Hoefler10,Peise12,Reisert17,Fukaya15, Fukaya18,Calotoiu2020}.
Our previous paper \cite{EigenKernel} proposed
several generic performance models 
in which the elapsed time $T$ of a routine is decomposed into the following five terms
\begin{eqnarray}
T_1(P)  &\equiv& \frac{c_1}{P}, 
\label{EQ-PERF-MODEL-TERM1} \\
T_2(P)  &\equiv& c_2, 
\label{EQ-PERF-MODEL-TERM2} \\
T_3(P)  &\equiv& c_3 \log P, \\
\label{EQ-PERF-MODEL-TERM3} 
T_4(P)  &\equiv& c_4 \frac{\log P}{\sqrt{P}},\  \\
\label{EQ-PERF-MODEL-TERM4} 
T_5(P)  &\equiv& \frac{c_5}{P^2}, 
\label{EQ-PERF-MODEL-TERM5} 
\end{eqnarray}
where the parameters $\{ c_i \}_{i}$  are to be estimated  
under the non-negative constraint ($\{ c_i \ge 0 \}_{i}$),
since each term is part of the elapsed time ($\{ T_i \ge 0 \}_{i}$). 
The terms $T_1$ and $T_2$ represent
the times of ideally parallel and non-parallel computations, respectively.
The two-term model $T =T_1 + T_2$ 
is the famous \rr{Amdahl's law} \cite{AMDAHL}.
The  term $T_3$ represents the \rr{setup} time of MPI communications
\cite{EigenKernel,PACHECO}. 
The three-term model $T = \sum_i^3 T_i$ is  
a minimal model that has its minimum at $P=P^{\ast} \equiv {\rm argmin}_P T(P)$.
The term $T_4(\propto \log P / \sqrt{P})$ expresses the time of MPI communications for matrix computation
\cite{EigenKernel} and 
the term  $T_5(\propto P^{-2})$ is \rr{a phenomenological term introduced to express} the \lq super-linear' behavior 
in which the time decays faster than $T_1 (\propto P^{-1})$ \cite{EigenKernel}.

This paper proposes 
an additional term responsible for the deceleration effect, 
when the total number of processor cores  exceeds the matrix size $M$. 
If the CPU  has $n_{\rm core}$ cores per node, 
the deceleration effect should appear for $P > P_{\rm c} \equiv M/n_{\rm core}$. 
The function for the deceleration effect can be given as
\begin{eqnarray}
T_6(P)  &\equiv& c_6 \frac{P}{1+ \exp \left( - (P-P_{\rm c}) \right)}, 
\label{EQ-PERF-MODEL-TERM6B} 
\end{eqnarray}
where a sigmoid-type function $(1+ \exp \left( - (P-P_{\rm c}) \right))^{-1}$ 
used as a \lq smoothed' step function. 
The term $T_6(P)$ is responsible for the deceleration effect
since the term vanishes ($T_6(P) \approx 0$) for $P \ll P_{\rm c}$ and 
is an increasing function ($T_6(P) \approx c_6 P$) for $P \gg P_{\rm c}$.

\subsection{Bayesian inference}

The parameter set $\{ c_i \}_{i=1,\nu}$ in the performance models $T=\sum_i^\nu T_i$ 
is estimated by Bayesian inference using the Monte Carlo (MC) method.
The posterior probability distribution is given by 
$\pi(X|D) \propto \pi(D|X) \pi(X)$,
where $D$ represents the teacher data set 
$\{ P_j, T_j^{\rm (exp)} \}_{j}$ and
$X \equiv (c_1, c_2, ..., c_\nu )$ is the parameter set to be estimated.
The cost function $F$ is defined from the relative error as 
\begin{eqnarray}
F \equiv \sum_j \frac{|T(P_j)-T_j^{\rm (exp)}|^2}{|T_j^{\rm (exp)}|^2}
\label{EQ-COST-FUNC}
\end{eqnarray}
and the likelihood $\pi(D|X)$ is defined to be 
proportional to $ \exp (- F/\tau)$
$(\pi(D|X) \propto  \exp (- F/\tau) )$.  
The given parameter $\tau(>0)$ is a measure of the tolerable uncertainty.
The non-negative constraint $(c_i \ge 0 )$ is reduced to
the uniform prior distribution in the region $[0, c_i^{\rm (max)}]$. 
The upper limit $c_i^{\rm (max)}$ is chosen to be sufficiently large so that the region of non-zero posterior distribution $(\pi(X|D)>0)$ 
is \rrr{contained} in the region $[0, c_i^{\rm (max)}]$. 
The cost function proposed here differs from 
that \rrr{used} in Ref.~\cite{EigenKernel}, as discussed in Sec. \ref{SEC-COST-FUNC}.

\subsection{Target problem and technical details}

We chose as our target problem  
the real-symmetric generalized eigenvalue problem 
using matrix data VCNT22500 from the ELSES matrix library
\cite{ELSES_MAT_LIB_URL,ELSES_MAT_LIB_PAPER},
where the matrix size $M$ is $M=22500$.  
The matrix data stem from the electronic state calculation by the ELSES simulator \cite{ELSES_URL, ELSES_PAPER}
of a vibrating carbon nanotube.  
The elapsed time data appear in the \lq total' time of Table 2 in Ref.~\cite{EigenKernel}, \rrr{which was measured on the K computer using ScaLAPACK}.
The data are available for various  numbers of nodes $P$, including  
$P=4, 16, 64, 256, 1024, 4096$ and 10,000.  
This data set was chosen 
because the elapsed time $T(P)$ has a minimum at $P=1024$,
which is typical in parallel computation. 
Since $n_{\rm core}=8$ and $M=22500$ in the present case, 
the critical node number $P_{\rm c}$ is $P_{\rm c} = 22500/8 = 2812.5$.

The MC method was carried out 
with the replica exchange Monte Carlo \rr{(REMC)} algorithm \cite{HUKUSHIMA-1996-REMC}
implemented in 
software package 2DMAT\cite{2DMAT-URL, 2DMAT-PAPER, HOSHI_2021_SION}.
The method is standard and the use of 2DMAT is not essential. 
The REMC method uses multiple values of the tolerable uncertainty parameters $\tau$, 
the chosen values of which are 
$\tau = \tau_0 \equiv 0.1 , \tau_1 \equiv 10^{2/3}\tau_0=0.464, \tau_2 \equiv 10^{4/3}\tau_0=2.154, \tau_3 \equiv 10^2\tau_0=10$. The total number of MCMC steps is $N_{\rm MCMC}=10^6$ \rrr{for each} value of the parameter $\tau$. 
The resultant sampling points with $\tau=0.1$ are used for the posterior probability density (histogram)
\rr{except for} the burn-in data of the early $N_{\rm MCMC}/2$ steps. 
The elapsed time for 2DMAT is approximately four minutes using a notebook computer.  

\begin{figure}[htb]
\centering
\includegraphics[width=7cm]{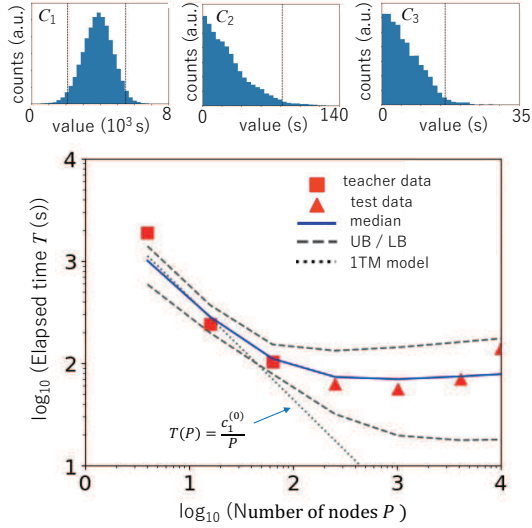}
\caption{Bayesian inference using the three-term model with three teacher data cases.
Upper panels: 
The posterior probability distribution (histogram) for parameters $c_1$ (upper left panel),
$c_2$ (upper middle panel), and $c_3$ (upper right panel). 
The upper and lower bounds of the 95 \% Highest Density Region (HDR) are drawn as vertical dashed lines in the upper left panel and
the upper bounds of the HDR are drawn in the upper middle and right panels. 
Lower panel: The elapsed times $T(P)$ for experimental data 
and the Bayesian inference.
The experimental data are shown for the three teacher data cases (squares) and the four test data cases (triangles). 
The results of 
the Bayesian inference are indicated
by the median (solid line) and the upper and lower bounds (UB/LB) of the 95 \% HDR (dashed lines). 
The time by the roughly estimated one-term (1TM) model ($T=c_{1}^{(0)}/P$) is  also shown (dotted line). 
\label{FIG-THREE-TERM-MODEL}
 }
\end{figure}

\section{Results \label{SEC-RESULTS}}

\subsection{Results with different models and teacher data \rr{sets}}

Figure \ref{FIG-THREE-TERM-MODEL} shows 
the predictions by the three-term (minimal) model ($T=\sum_i^3 T_i$)
with teacher data cases at $P=4, 16, 64$. 
The remainder of the experimental data, at $P=256, 1024, 4096, 10000$, was \rr{used as} the test data cases. 
The posterior probability density
is shown in the histograms for $\{ c_i \}_{i}=1,2,3$ in the left, middle and right panels, respectively, 
of the upper section of Fig.~\ref{FIG-THREE-TERM-MODEL}. 
The histograms for $c_1$, $c_2$ and $c_3$ show the maximums located near $c_1 \approx  c_1^{\rm (0)} \equiv 4 \times 10^3$ s, $c_2 \approx 0$ and $c_3  \approx 0$, respectively, which leads us to the one-term model $T \approx c_1^{\rm (0)} / P$ as a rough estimation. 
The maximums for $c_2$ and $c_3$ are located near the origin, owing to the non-negative constraint ($\{ c_i \ge 0 \}_i$).
In the lower panel of Fig. ~\ref{FIG-THREE-TERM-MODEL}, 
the elapsed times $T(P)$ \rr{of}
the experimental data \rr{are compared with} those \rr{by} the Bayesian inference.
The experimental data are \rr{contained} in the 95 \% HDR \rr{(Highest Density Region)}, except at $P=4$. 
The Bayesian inference indicates that parallelism with $P \ge 1000$ is not an efficient parallel computation,
which offers a fruitful guideline for application researchers.

\begin{figure*}[ht]
\centering
\includegraphics[width=16cm]{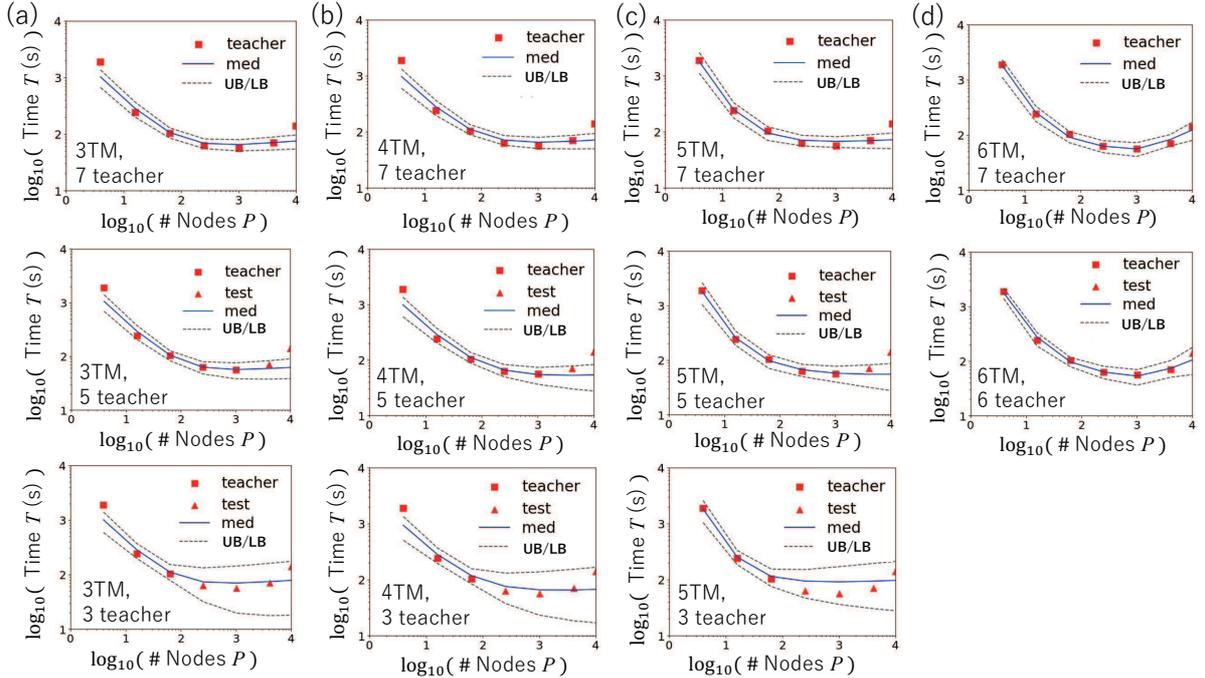}
\caption{Comparison of the elapsed times $T(P)$ for the experimental data 
and for the Bayesian inference. 
The experimental data are shown 
as squares for the teacher data cases and as triangles for the test data cases.   
The results of the Bayesian inference 
are indicated 
by the median (solid line) and the upper and lower bounds (UB/LB) of the 95 \% Highest Density Region (dashed lines). 
The Bayesian inference was carried out
using 
(a) the three-term (3TM) model with  the seven, five and three teacher data cases, 
(b) the four-term (4TM) model with  the seven, five and three teacher data cases,  
(c) the five-term  (5TM) model with  the seven, five and three teacher data cases, and
(d) the six-term  (6TM) model with  the seven and six teacher data cases. 
\label{FIG-COMPARISON}
}
\end{figure*}


Figure ~\ref{FIG-COMPARISON} shows a comparison 
of the  Bayesian inference with different models and teacher data cases.
Here, it is evident  that the width of the 95 \% HDR in the test data region is wider than in the teacher data region,
as would be expected.  
In the lower panel of Fig.~\ref{FIG-COMPARISON}(a), 
for example, the widths of the 95 \% HDR at $P=256, 1024, 4096, 10000$ are wider than those at $P=4,16,64$.
Figures ~\ref{FIG-COMPARISON}(a), (b) and (c) show the results of the three-term model ($T=\sum_i^3 T_i$),
the four-term model ($T=\sum_i^4 T_i$) and the five-term model ($T=\sum_i^5 T_i$), respectively.
The upper, middle and lower panels show the results with seven, five and three teacher data cases, respectively.
Figure ~\ref{FIG-COMPARISON}(d) shows the results of the six-term model 
with seven and six teacher data cases
in the upper and lower panels, respectively.
Since $T_6(P) \approx 0$ for $P \ll P_{\rm c}=2812.5$, the six-term model is meaningful only when the data case at $P=4096$ and/or $P=10000$ are included in the teacher data. 

Figures \ref{FIG-COMPARISON}(a) and (b) indicate that the median curves in the three- and four-term models
reproduce the experimental data \rr{satisfactorily}, except at $P=4$ and $P=10000$. 
The deviations at $P=4$ and $P=10000$ stem from the limited representation ability of these models,
since the deviations appear even with the seven teacher data cases in the upper panels of Figures ~\ref{FIG-COMPARISON}(a) and (b). 
The \rr{deviation} at $P=4$ disappears with the five-term model, as shown Figure ~\ref{FIG-COMPARISON}(c),
since the rapid decrease from $T(P=4)$ \rrr{to} $T(P=16)$ can be expressed by the super-linear term ($T_5$).
The \rr{deviation} at $P=10000$ disappears with the six-term model and the six teacher data cases, as shown in Figure ~\ref{FIG-COMPARISON}(d),
since the rapid increase from $T(P=4096)$ into $T(P=10000)$ can be expressed by the deceleration term ($T_6$).


\subsection{Comparison \rr{of} different cost functions \label{SEC-COST-FUNC}}

We can now make a comparison \rr{of} the different cost functions. 
In the previous paper \cite{EigenKernel}, 
transformed variables $(\rho, \phi)  \equiv (\log P, \log T)$ are used 
rather than the original variables $(P, T)$ 
and the cost function is defined as the absolute error of 
$F_{\rm trf} \equiv \sum_j |\phi(\rho_j)-\phi_j^{\rm (exp)}|^2$,  unlike in Eq.(\ref{EQ-COST-FUNC}).

It is possible to compare the results produced by the different cost functions,
since the results from the three-term model with three teacher data cases are given  
in the lower panel of Fig.~\ref{FIG-COMPARISON}(a) and Fig. 6(a) of Ref.~\cite{EigenKernel}
and the results from the five-term model with three teacher data cases are given
in the lower panel of Fig.~\ref{FIG-COMPARISON}(c) and Fig. 6(b) of Ref.~\cite{EigenKernel}.
The comparison \rr{suggests} that 
the results using the present cost function $F$ are significantly better 
than those with the previous cost function $F_{\rm trf}$.

\subsection{Discussion}

It should be recalled that general  
application researchers are interested primarily in the location of the minimum point ($P=P^\ast \equiv {\rm argmin}_P T(P) \approx 10^3$). 
In addition, application researchers have a solid need 
for a simple performance model in which  the number of parameters is small 
and the parameters can be estimated from a small teacher data set without overfitting. 
Figure ~\ref{FIG-COMPARISON} indicates 
that 
the three-term model with the three teacher data cases \rrr{gives satisfactory prediction} at $P \approx 10^3$, which leads to the conclusion that  
the minimal three-term model with  three teacher data cases is effecctive.  

A possible future direction for model development is prediction (extrapolation) with respect to \rrr{both}
the number of nodes $P$ and the matrix size $M$ ($T=T(P,M)$).
The modeling will be fully realized,
when the the coefficients $\{ c_i \}$ in 
the term $\{ T_i \}$
are re-defined as a function of $M$ ($c_i \equiv c_i(M)$).
A reasonable modeling for matrix computation is a third order polynomial
$(c_i(M) \equiv c_{i3} M^3 + c_{i2} M^2 +c_{i1} M +c_{i0})$. 
In such a case, 
the number of the parameters ($\{ c_{ij} \}_{ij}$)
in the five-term model ($4 \times 5 = 20$) 
is significantly larger than in the three-term model ($4 \times 3 = 12$). 
The above discussion implies the importance of the three-term (simpler) model. 

It would also be desirable to \rrr{obtain more detailed performance data by measuring} the elapsed time \rrr{of lower-level routines}.
For example, one could measure the execution time of \rr{an} MPI \rr{function},
by \rr{inserting} the timer \rr{into} the source code.
Such a detailed elapsed time could be \rr{used to construct} a \rrr{more accurate and reliable} performance model for the \rrr{entire solver}.

\section{Summary \label{SEC-SUMMARY}}

A performance prediction method for parallel computations is proposed
using \rr{parametrized} performance \rr{models} and Bayesian inference.
The proposed method meets the need of application researchers,
in particular, for performance extrapolation when the number of used nodes is large. 
Three-, four-, five-, and six-term performance models were proposed and applied to the real-symmetric generalized eigen-value problem.
Comparison of the models indicates that the proposed method is effective even with the \rrr{the combination of} the three-term (minimal) model and the minimal teacher data set. 
Possible directions for further development of the method were suggested,
including extrapolation 
with respect to both the number of nodes $P$ and the matrix size $M$.
Importantly, the prediction method is general and can be used in 
applications other than the generalized eigenvalue problem.

\acknowledgments

The present research is supported by the Grant-in-Aid
for Scientific Research (KAKENHI) from Japan Society for the
Promotion of Science (19H04125, 20H00581) and JHPCN and HPCI in Japan (jh210044-NAH).

\references



%
%
%
%

\end{document}